\documentclass[A4,12pt]{article}

\usepackage{hyperref}

\usepackage{latexsym}
\usepackage{amssymb,amsmath}

\newcommand{\be}{\begin{equation}}      
\newcommand{\ee}{\end{equation}}        
\newcommand{\bern}{\begin{eqnarray*}}   
\newcommand{\eern}{\end{eqnarray*}}     
\newcommand{\beqp}{\begin{eqproof}}     
\newcommand{\eeqp}{\end{eqproof}}       

\newcommand{\bt}{\begin{teorema}}       
\newcommand{\et}{\end{teorema}}         
\newcommand{\bd}{\begin{definizione}}   
\newcommand{\ed}{\end{definizione}}     
\newcommand{\bc}{\begin{corollario}}    
\newcommand{\ec}{\end{corollario}}      
\newcommand{\bp}[1]{\noindent{\bf Proof #1.} } 
\newcommand{\ep}{\hfill$\Box$ }   

\newcommand{\bl}{\begin{lemma}}         
\newcommand{\el}{\end{lemma}}           
\newcommand{\boss}{\begin{osservazione}\rm}     
\newcommand{\eoss}{\end{osservazione}}  

\font\corsivo=rsfs10 at 12pt
\font\doppio=msbm10 at 12pt
\font\scdoppio=msbm8

\newcommand{\C}{\mbox{\corsivo C}}      
\newcommand{\R}{\hbox{\doppio R}}       
\newcommand{\RR}{\hbox{\scdoppio R}}    
\newcommand{\G}{\hbox{\doppio G}}     

\newcommand{\abs}[1]{\left|#1\right|}   
\newcommand{\decl}{:=}                  

\newtheorem{teorema}{Theorem}[section]
\newtheorem{definizione}[teorema]{Definition}
\newtheorem{corollario}[teorema]{Corollary}
\newtheorem{lemma}[teorema]{Lemma}

\newtheorem{example}[teorema]{Example}

\newtheorem{osservazione}[teorema]{Remark}

\def\acknowledgment{\goodbreak\subsection*{Acknowledgment}
\bgroup \footnotesize}

\def\acknowledgments{\goodbreak\subsection*{Acknowledgments}
\bgroup \footnotesize}

\def\endacknowledgment{\vskip1sp\egroup}
\def\endacknowledgments{\vskip1sp\egroup}

\def\fix{\ifcase\oldtime 0\or0\or0\or%
0\or0\or0\or0\or0\or0\or0\fi}
\def\fixtiming{\ifcase\timing 0\or0\or0\or%
0\or0\or0\or0\or0\or0\or0\fi}

\newcount\timing
\newcount\hourcount
\newcount\oldtime

\def\realtime{\timing=\time \oldtime=\time
\ifnum\timing>60 \divide\timing by 60
\hourcount=\the\timing
\multiply\timing by 60
\advance\oldtime by-\timing
\ifnum\hourcount<12 \number\hourcount:\fix\number\oldtime am\fi%
\ifnum\hourcount=12 \number\hourcount:\fix\number\oldtime pm\fi%
\ifnum\hourcount>12 \advance\hourcount by-12
\number\hourcount:\fix\number\oldtime pm\fi
\else12:\fixtiming\number\timing am\fi}

\newcommand{\diver}{\mathrm{div}}       

\newcommand{\Cuno}{\mbox{\corsivo C}^{\,1}}

\newcommand{\hei}{\hbox{\doppio H}}     

\newcommand{\RN}{\R^N}

\newcommand{\grl}{\nabla_{\!\!L}}               
\newcommand{\diverl}{\diver_{\!L}}      

\newcommand{\grh}{\nabla_{\!\!H}}               
\newcommand{\lh}{\Delta_H}              

\newcommand{\email}[1]{\tt #1}

\begin{document}
\date{May 7 2009\\ \small Accepted for publication on \\
Mat. Sb., 2010, Vol. 201 }
\title{Nonnegative solutions of some \\quasilinear elliptic inequalities
and applications}
\author{Lorenzo D'Ambrosio\footnote{Dipartimento di Matematica,
via E. Orabona, 4 - 
Universit\`a degli Studi di Bari,
I-70125 Bari,
Italy
\email{dambros@dm.uniba.it}} \ \ and\ \  
Enzo Mitidieri\footnote{Dipartimento di Matematica e Informatica,
via A. Valerio, 12/1 - 
Universit\`a degli Studi di Trieste,
I-34127 Trieste,
Italy
\email{mitidier@units.it}}}

\maketitle
\begin{abstract}
Let $f :\R\to\R$ be a continuous function. We prove
that under some additional assumptions on $f$ and $A:\R\to\R_{+}$, weak $\Cuno$ solutions of the 
differential inequality
$-\diver(A(\abs{\nabla u})\nabla u)\ge f(u)$ on $\RN$
are nonnegative.
Some extensions of the result in the framework of subelliptic
operators on Carnot Groups are considered.
\end{abstract}

\section{Introduction}\label{sec:intrhi}

In this paper we shall study the following problem.  

Let $L$ be a second order differential operator and let $f:\R\to\R$ be a continuous function.
Find  additional assumptions on $(L,f)$ that imply the positivity of the possible solutions of the
differential inequality
\be
L(u)\ge f(u) \quad\mathrm{on}\quad\RN.\label{cond:0}
\ee

Some partial answers to this problem have been obtained  in \cite{dam-mit-poh04, mit-poh03}. In those papers, the authors deal with elliptic inequalities of the  form (\ref{cond:0}) in the case when $L$ is 
the Laplacian operator or the polyharmonic operator $(-\Delta)^k$ in the Euclidean setting or, more generally $L$ is a sub elliptic Laplacian on a Carnot group and $f$ is nonnegative.
The main strategy used in \cite{dam-mit-poh04, mit-poh03} for proving positivity results was
via integral representation formulae. One essential difficulty using this approach is that no assumptions  on the behavior of the solutions at infinity are known. A typical example in this direction is given by,
\be-\Delta u \ge \abs u ^q\quad\mathrm{on}\quad\RN,\label{cond:1}
\ee
where $N\ge3$ and $q>1$.

\medskip
The following result holds (see \cite{dam-mit-poh04}).

\bt\label{th:model} Let $N\ge3$ and $ q>1$. Let $u\in L^q_{loc}(\RN)$ be a distributional solution of (\ref{cond:1}) and let $Leb(u)$ be the set of its Lebesgue points.
If $x\in Leb(u)$, then 

$$u(x)\ge C_{N} \int_{\RR^{N}}\frac{\abs {u(y)}^q}{|x-y|^{N-2}}\,\,\,dy,$$
where $C_{N}$ is an explicit positive constant.

\et

From this result it follows that, if $u$ is a solution of (\ref{cond:1}) then, either $u(x)=0$ 
or $u(x)>0$ a.e. on $\RN.$

Obviously, the approach via representation formulae  cannot be applied to quasilinear problems. 
In this paper we shall consider a class of quasilinear model problems for which the positivity property
mentioned at the beginning of this introduction holds. 

More precisely, we shall deal with the case when $L$ is the $p$-Laplacian operator,
namely $\Delta_p\cdot=\diver(\abs{\nabla \cdot}^{p-2}\nabla \cdot)$, or the mean curvature operator
$\displaystyle -\diver(\frac{\nabla \cdot}{\sqrt{1+\abs{\nabla \cdot}^2}}) $.
In this cases some results on nonnegativity of solutions of (\ref{cond:0}) are proved by using a suitable comparison Lemma (see Lemma \ref{lem:comp} below).
We will apply this results to the problem of  a priori bounds of solutions and nonexistence theorems. For  interesting results on related coercive equations see the very recent paper \cite{FAS}.

This paper is organized as follows. In section $2$ we state and prove our main results in the Euclidean setting
for the $p$-Laplacian and the mean curvature operator.
In addition we  point out some consequences and briefly discuss the sharpness of some of the results.

In section $3$ we briefly indicate some generalizations to other quasilinear operators
and to differential inequalities on Carnot groups. We end the paper with an appendix which contains some well known facts on Carnot groups.

A shorter version of  this paper  appeared in \cite{dammit09}.
\noindent
After this paper was submitted for publication (May 2009), we learned form Professor James Serrin that in the Euclidean setting, related results similar
to Lemma \ref{lem:comp}  are contained in \cite{PS}. See in particular  Sections 2.2 and 2.4 and also Chapter 3. It seems likely that
suitable versions of the comparison principles contained in \cite{PS} still hold in the framework of Carnot groups. 
This will be investigated in a forthcoming paper.
\section{Main results}\label{sec:main}
Throughout this  paper we shall assume that $N\ge 2 $ and  we deal with weak $\Cuno$ solutions of the problems under consideration.
 See Definition \ref{def:sol} below for further details.
Our main results are the following.

\bt\label{th:main} Let $p>1$ and $N>1$. Let $f\colon\R\to\R$ be a continuous function such that 
\be f(t)>0\quad \mathrm{ if}\quad  t<0, 
\quad f\quad \mathrm{is\ non\ increasing\ on}\ \ ]-\infty, 0 [\label{cond:f}\ee 
and
  \be \int_{-\infty}^{-1} \left(\int_t^{-1} f(s)\ ds\right)^{-\frac1p} \ dt< +\infty. \label{dis:naiusa}\ee
  If $u$ is a solution of 
  \be -\diver(\abs{\nabla u}^{p-2} \nabla u)\ge f(u)\quad\mathrm{on} \quad \RN, \label{dis:main}\ee
  then $u\ge0$ on $\RN$. Moreover if $f(t)\ge 0$ for $t\ge 0$ then, either $u\equiv 0$ or $u>0$ on $\RN$.
\et 

\bc  Let $p>1$. Let $f\colon\R\to\R$ be a continuous function such that 
  $f(t)\ge C \abs t^q$ for $t<0$.
  Let $u$ be a solution of 
 \be -\diver(\abs{\nabla u}^{p-2} \nabla u)\ge f(u)\quad\mathrm{on} \quad \RN. \label{dis:mainpower}\ee
 If $q>p-1$ then $u\ge0$ on $\RN$. Moreover if $f(t)\ge 0$ for $t\ge 0$ then, either $u\equiv 0$ or $u>0$ on $\RN$.
\ec

In the case of the mean curvature operator the above results can be improved. 
Indeed, the claim follows without the assumption
(\ref{dis:naiusa}) on $f$.
\bt\label{th:euclmean} Let  $f\colon\R\to\R$ be a continuous function 
 satisfying (\ref{cond:f}). Let u be a solution of
\be -\diver(\frac{\nabla u}{\sqrt{1+\abs{\nabla u}^2}})\ge f(u)\quad\mathrm{on} \quad \RN. \label{dis:mainmean}\ee
Then $u\ge 0$ on $\RN$.
\et

A first consequence of the above results  is the following a priori estimate.
\bt\label{th:eqpbd}  Let $p>1$ and $N>1$. Let $f\colon\R\to\R$ be a continuous function such that there exists $\alpha,\beta\in\R$, $\alpha\le\beta$ such that
\be f_{]-\infty,\alpha[}\quad \mathrm{ is\ positive\ and\ non\ increasing,}\quad 
f_{]\beta,+\infty[}\quad \mathrm{ is\ negative\ and\ non\ increasing}, 
    \label{cond:fodd}\ee
and
\be \int_{-\infty}^{\alpha} \left(\int_t^{\alpha} f(s)\ ds\right)^{-\frac1p} \ dt< +\infty , \quad
   \int_\beta^{\infty} \left(\int_\beta^{t} -f(s)\ ds\right)^{-\frac1p} \ dt< +\infty. \label{dis:naiusa00}\ee
  If $u$ is a solution of 
  \be -\diver(\abs{\nabla u}^{p-2} \nabla u)= f(u)\quad\mathrm{on} \quad \RN, 
    \label{eq:pgen}\ee
  then $u$ is bounded and  $\alpha\le{u(x)}\le \beta$  for any $x\in \RN$.
\et

Again, for the mean curvature operator we can require more general assumption on $f$. 
\bt\label{th:eqmeanbd}  
  Let $f\colon\R\to\R$ be a continuous function such that 
  $$\liminf_{t\to-\infty} f(t)>0.$$
  If $u$ is a solution of (\ref{dis:mainmean}), then
  $f$ has at least a zero, and set $\alpha:=$ the first zero of $f$
  (that is $\alpha:=\min S$ where $S:= f^{-1}(0)$) we have $u\ge \alpha$.
  In particular if $f>0$ the (\ref{dis:mainmean}) has no solution.

  Moreover if   $$\limsup_{t\to+\infty} f(t)<0$$
  and $u$ solves
 \be -\diver(\frac{\nabla u}{\sqrt{1+\abs{\nabla u}^2}})= f(u)\quad\mathrm{on} \quad \RN. \label{eq:mean1}\ee
  then $u$ is bounded and  $\alpha\le{u(x)}\le \beta$  for any $x\in \RN$,
  where $\beta:=$ last zero of $f$ (that is $\beta:=\max S$).
\et

A direct consequence of  Theorems \ref{th:eqpbd} and \ref{th:euclmean}, \ref{th:eqmeanbd} 
are the following Liouville theorems.
\bc\label{th:eqp}  Let $p>1$ and $N>1$. Let $f\colon\R\to\R$ be a non increasing continuous function such that 
\be f(t)>0\quad \mathrm{ if}\quad  t<0, \ \  \mathrm{and}\  \  f(t)<0\quad \mathrm{ if}\quad  t>0, \label{cond:odd}\ee
and
\be \int_{-\infty}^{-1} \left(\int_t^{-1} f(s)\ ds\right)^{-\frac1p} \ dt< +\infty , \quad
   \int_1^{\infty} \left(\int_1^{t} -f(s)\ ds\right)^{-\frac1p} \ dt< +\infty. \label{dis:naiusa2}\ee
  If $u$ is a solution of 
  \be -\diver(\abs{\nabla u}^{p-2} \nabla u)= f(u)\quad\mathrm{on} \quad \RN, 
    \label{eq:plapg}\ee
  then $u\equiv 0$ on $\RN$.

In particular, if $q>p-1$ and $u$ is a solution of 
  \be \Delta_p u= \abs{u}^{q-1}u\quad\mathrm{on} \quad \RN, \label{brezisp}\ee
  then $u\equiv 0$ on $\RN$.
\ec

\boss 
The conclusion for equation (\ref{brezisp}) in the case $p=2$ has been proved by Brezis \cite{brezis}.
An important generalization, by using a different technique,
in the Euclidean case,
has been obtained by Serrin \cite[Theorem 2]{serrin1}. For a version of Corollary \ref{th:eqp} in the framework of Carnot groups
see Theorem \ref{th:eqpG} below.
\eoss

\bt\label{th:disp} Let $p>1$ and $N>1$. Let $f\colon\R\to\R$ be a positive, non increasing, continuous function 
satisfying (\ref{dis:naiusa}).
Then the inequality (\ref{dis:main}) has no solutions.
\et

\bt\label{th:meaneq}  Let $N>1$. Let $f\colon\R\to\R$ be a non increasing,
  continuous function and $f\not\equiv 0$. If $u\in C^2( \RN)$ is a solution of 
\be -\diver(\frac{\nabla u}{\sqrt{1+\abs{\nabla u}^2}})= f(u)\quad\mathrm{on} \quad \RN \label{eq:mean2}\ee
then $u$ is constant, that is $u\equiv \alpha$ and $f(\alpha)=0$.
In particular if $f(t)\neq 0$ for any $t$, then (\ref{eq:mean2}) has no solutions.

In addition, if $f$ is supposed to be positive then the inequality
\be -\diver(\frac{\nabla u}{\sqrt{1+\abs{\nabla u}^2}})\ge f(u)\quad\mathrm{on} \quad \RN\label{dis:mean33} \ee
has no solutions.
\et
\boss For a different proof of the first part of  Theorem \ref{th:meaneq} 
 see  \cite{FA}. 
  The claim concerning the inequality (\ref{dis:mean33})
  is new and of independent interest. 
\eoss

The following results are an easy consequence of the fact that the only
nonnegative functions $u$ such that $-\Delta_p u\ge0$ on $\RN$ with $N\le p$ or 
 $\displaystyle -\diver(\frac{\nabla u}{\sqrt{1+\abs{\nabla u}^2}})\ge0$ on $\R^2$ are the 
constants, see \cite{MP2}.
\bc\label{cor:liouv} Let $p\ge N>1$ and $f:\R\to[0,+\infty[$ be  a continuous function satisfying (\ref{cond:f}) and
(\ref{dis:naiusa}). If $u$ is a solution of  (\ref{dis:main}) then $u$ is constant on $\RN$.
More precisely $u\equiv \alpha\ge 0$ and $f(\alpha)=0$.

Moreover  if $f(t)>0$ for $t\ge0$, then  (\ref{dis:main}) has no solutions.
\ec
\bc Let $f:\R\to[0,+\infty[$ be  a continuous function satisfying (\ref{cond:f}). 
  Let u be a solution of
\be -\diver(\frac{\nabla u}{\sqrt{1+\abs{\nabla u}^2}})\ge f(u)\quad\mathrm{on} \quad \R^2. \label{dis:meanr2}\ee 
Then $u$ is constant on $\R^2$. More precisely 
$u\equiv \alpha\ge 0$ and $f(\alpha)=0$.

Moreover if $f(t)>0$ for $t\ge0$, then  (\ref{dis:meanr2}) has no solutions.
\ec

\boss The above assumptions on $f$ are sharp in the following sense. If $p=2$ and $q=1=p-1$
  the result is false. Indeed the equation
  $$-\Delta u =\abs u\quad\mathrm{on}\quad\RN, $$ admits the explicit negative  solution
  $$u(x):= -Exp(x_{1}),\quad x\in \RN, $$ or solutions that changes sign (see \cite{dam-mit-poh04}).

  In the general case $q=p-1$ the equation
  $$\Delta_p u= u^{p-1}\quad\mathrm{on}\quad\RN,$$
  admits a positive  solution (see for instance \cite{nai-usa97}). Therefore,
  the equation  $$-\Delta_p u=\abs u^{p-1} \quad\mathrm{on}\quad\RN, $$
  has a negative solution.
\eoss

Let us briefly describe the idea of the proof of our main result.
Let $u$ be a solution of (\ref{dis:main}). Without loss of generality we will
show that $u(0)\ge 0$. The function $U\decl-u$ satisfies the inequality 
$$\diver(\abs{\nabla U}^{p-2} \nabla U)\ge f(-U)\quad\mathrm{on} \quad \RN.$$ 
Let $v$ be a positive solution of 
$$\diver(\abs{\nabla v}^{p-2} \nabla v)= f(-v)\quad\mathrm{on} \quad B_R,$$ 
such that $v(0)=a>0$ and $v(x)\to+\infty$ as $\abs x\to R$.
The assumptions on $f$ imply the existence of $v$.
Since $U(x)\le v(x)$ for $\abs x$ close to $R$, by a comparison Lemma
(see Lemma below) it follows that $U(x)\le v(x)$ for any $\abs x< R$. In particular $U(0)\le v(0)=a$.
Letting $a\to 0$ we have $U(0)\le 0$. Hence $u(0)\ge 0$.
Finally, if $f(t)\ge0$ for $t\ge0$, by the weak Harnack inequality we get that, either $u\equiv 0$ or $u> 0$ on $\RN$.

\subsection{A comparison Lemma}

In this section, we shall prove a comparison Lemma. This Lemma will be useful when dealing
with more general operators then the $p$-Laplacian and the mean curvature operator, 
thus we shall present it in a general form.
To this end let us introduce some notations. Let $\mu=(\mu(x)_{i,j})$ be a matrix 
with $N$ columns and $l(\le N)$ rows with entries belonging, for simplicity, to $\Cuno(\RN)$.
We denote by $\grl\decl\mu\nabla$ and by $\diverl=-\grl^*=\diver(\mu^T\cdot)$.
The isotropic Euclidean case corresponds to the choice $\mu=I_N$ the unitary matrix of dimension $N$.

We shall assume that if $\grl u=0$ on a connected region $\Omega$ then $u\equiv const$ in such a region.

\bd\label{def:sol} Let $\Omega\subset\RN$ be an open set, let $A:\R\to\R$ and $h:\Omega\to\R$ be 
 a continuous functions. We say that $u$ \emph{is a solution} of 
$$ \diverl(A(\abs{\grl u})\grl u)\ge h \quad\mathrm{on}\quad\Omega,$$
if $u\in\Cuno(\Omega)$  and for 
any nonnegative $\phi\in\Cuno_0(\Omega),$ we have
$$ -\int_\Omega A(\abs{\grl u})\grl u \cdot\grl\phi \ge\int_\Omega h\phi. $$

In a similar manner we can define  solutions of the inequalities $$- \diverl(A(\abs{\grl u})\nabla u)\ge h\ \  and\ \  
\diverl(A(\abs{\grl u})\nabla u)\le h.$$
\ed

The following Lemma  is useful when considering
solutions of inequalities of the form,
\be \diverl(A(\abs{\grl u})\grl u)\ge g_1(x,u)\quad\mathrm{on}\quad\Omega,
                      \label{dis:coeu}\ee  
and
\be \diverl(A(\abs{\grl v})\grl v)\le g_2(x,v)\quad\mathrm{on}\quad\Omega.
                      \label{dis:coev}\ee  
Here,  $A$ is a continuous function such that $A(t)>0$ for $t>0$ and for  $i=1,2$,
$g_i:\Omega\times\R\to\R$ is continuous.

\bl\label{lem:comp} Let $\Omega$ be a bounded open set and let
   $u$ and $v$ be solutions of (\ref{dis:coeu}) and
  (\ref{dis:coev}) respectively. Assume that
  \begin{enumerate}
  \item
    \begin{enumerate}
    \item For any $x\in\Omega$, $t\ge s\ge0$ there holds $g_1(x,t)\ge g_2(x,s)$,
      $g_1(x,\cdot)$ is not decreasing on $]0,+\infty[$ and $v\ge 0$;
    \item[or]
      \item For any $x\in\Omega$, $t\ge s $ there holds $g_1(x,t)\ge g_2(x,s)$ and
      $g_1(x,\cdot)$ is not decreasing;
    \end{enumerate}
  \item The function $tA(t)$ is increasing and positive for $t>0$;
  \item $u\le v$ on $\partial\Omega$.
  \end{enumerate}
 Then $u\le v$ on $\Omega$.
\el

\bp{} Let $u$ and $v$ be solutions of (\ref{dis:coeu}) and (\ref{dis:coev}) respectively. 
  Let $\epsilon>0$ be fixed and set $v_\epsilon\decl v+\epsilon$.
  It is a simple matter  to check that the function $v_\epsilon$ satisfies  the inequality
  $$ \diverl(A(\abs{\grl w})\grl w)\le g_1(x,w)\quad\mathrm{on}\quad\Omega. $$

  Therefore, for any nonnegative $\phi\in\Cuno_0(\Omega)$ we have

  \be-\int_\Omega(A(\abs{\grl u})\grl u -A(\abs{\grl v_\epsilon})\grl v_\epsilon)
    \cdot\grl\phi \ge\int_\Omega (g_1(x,u)-g_1(x,v_\epsilon))\phi. \label{dis:tec1}\ee
 Next we choose $\phi$ as follows: $\phi\decl ((u-v_\epsilon)_+ )^2$.
  It is clear that $\phi$ is nonnegative and $\phi\in\Cuno(\Omega)$. Moreover, since 
    $v_\epsilon-u\ge\epsilon>0$ on $\partial\Omega$, it follows that $\phi$ has compact support.
  Substituting $\phi$ in (\ref{dis:tec1}), we obtain
   $$  -\int_\Omega(A(\abs{\grl u})\grl u -A(\abs{\grl v})\grl v)
    \cdot(\grl u-\grl v) 2 (u-v_\epsilon)_+ 
    \ge\int_\Omega (g_1(x,u)-g_1(x,v_\epsilon))\phi.$$

  We claim that 
  \be (A(\abs{\grl u})\grl u -A(\abs{\grl v})\grl v)\cdot(\grl u-\grl v)\ge0. 
              \label{dis:claim1} \ee
  Indeed,
  \begin{eqnarray}
   \lefteqn{ (A(\abs{\grl u})\grl u -A(\abs{\grl v})\grl v)\cdot(\grl u-\grl v)}\nonumber \\
   &&=A(\abs{\grl u})\abs{\grl u}^2+A(\abs{\grl v})\abs{\grl v}^2- 
        (A(\abs{\grl u})+A(\abs{\grl v}))(\grl u\cdot\grl v)\nonumber \\
   &&= \Big(A(\abs{\grl u})\abs{\grl u}-A(\abs{\grl v})\abs{\grl v}\Big)
        \Big(\abs{\grl u}-\abs{\grl v}\Big)+ \nonumber\\
   &&\qquad + \Big( A(\abs{\grl u})+A(\abs{\grl v})\Big)\Big(\abs{\grl u}\abs{\grl v}-\grl u\cdot\grl v\Big)=:I_1+I_2.  \label{eq:cmpt}
    \end{eqnarray}
  Since $A\ge 0$, we have $I_2\ge 0$. 
  From the monotonicity of $tA(t)$, it follows that
  $$ I_1= \Big(A(\abs{\grl u})\abs{\grl u}-A(\abs{\grl v})\abs{\grl v}\Big)
        \Big(\abs{\grl u}-\abs{\grl v}\Big)
 \ge 0. $$

  Assume first that  $g_{1}(x,.)$ is strictly increasing.
  Therefore, 
  the inequality 
  $(g_1(x,u)-g_1(x,v_\epsilon) ((u-v_\epsilon)_+ )^2\ge0$
  holds for every $x\in\Omega$. As a consequence, 
  $(g_1(x,u)-g_1(x,v_\epsilon) ((u-v_\epsilon)_+ )^2=0$ on $\Omega$ and hence $u\le v_\epsilon$.
  
  This completes the proof in case $g_{1}(x,.)$ is strictly increasing.
For the general case we need of an extra argument. Indeed,
from (\ref{dis:claim1}) and (\ref{eq:cmpt}) we have that
 $\int_\Omega (I_1+I_2)(u-v_\epsilon)_+=0$.
  Let $x\in\Omega$ be such that $u(x)\ge v_\epsilon(x)$. Since $I_1\ge0$ and $I_2\ge 0$
  we have $I_1(x)=0=I_2(x)$. 

  We claim that $\grl u (x) =  \grl v(x)$. Indeed, if  $\grl u (x) \neq  \grl v(x)$,
  from $I_2(x)=0$, we deduce that 
  $\abs{\grl u (x)} \neq  \abs{\grl v(x)}$\footnote{If $t,s$ are two different vectors 
    in a Hilbert space such that $(s\cdot t)=\abs t \abs s$, then 
    $0<\abs{t-s}^2=\abs t^2 +\abs s^2-2(s\cdot t)=\abs t^2 +\abs s^2-2 \abs s \abs t=
    (\abs s -\abs t)^2 $.}.
  Thus from $I_1(x)=0$, the fact that $tA(t)$ is injective, it follows that
  $$0=A(\abs{\grl u})\abs{\grl u}-A(\abs{\grl v})\abs{\grl v} \neq 0.$$
This implies
  $\grl ((u-v_\epsilon)_+ )^2=0$ on $\Omega$, that is $((u-v_\epsilon)_+ )^2=\phi=0$.
  
  Therefore, letting $\epsilon\to 0$ in $u\le v+\epsilon$ the claim follows.
\ep
\boss It is possible to extend the above result to situations where the function $A$  depends on $x$.
  Namely, when $A:\Omega\times\R\to \R$ is continuous and for any $x\in \Omega$
  the function $t\in]0,+\infty[\to tA(x,t)$ is increasing and strictly positive.
\eoss
\subsection{Proofs of the Main Results}
The proof of Theorem \ref{th:main} relies on the following.

\bt\label{th:blowup} Let $g\colon\R\to\R$ be a continuous function such that, 
  $g(t)>0$ if  $t<0$, $g$ is non decreasing on $]0,+\infty[$, and
  \be \int_1^{+\infty} \left(\int_1^t g(s)\ ds\right)^{-\frac1p} \ dt< +\infty. \label{dis:ossp}\ee
  For any $p>1$, $a>0$, $D>1$, there exists a function $\varphi$ and $R>0$ such that, 
  $\varphi$ is a solution of
  \be \left(r^{D-1} \abs{\varphi'(r)}^{p-2}\varphi'(r)\right)'=r^{D-1} g(\varphi(r)),
      \quad \varphi(0)=a,\quad \varphi'(0)=0, \label{eq:radial} \ee 
  $\phi$ is increasing on $]0,R[$ and $\varphi(r)\to +\infty$ as $r\to R$.
\et
  See \cite{oss57} for a proof in the case $p=2$ and \cite{nai-usa97} for the quasilinear 
  case $p\neq 2$.

\bigskip
\bp {of Theorem \ref{th:main}}
Let $u$ be a solution of (\ref{dis:main}). 
Since the inequality is invariant under translations, it is sufficient to prove that $u(0)\ge 0$.

Let $g(t)\decl f(-t)$. The function $g$ satisfies the assumptions of Theorem \ref{th:blowup}.
Let $D=N>1$, $a>0$ and let $\varphi$ be a solution of (\ref{eq:radial}) such that 
$\varphi(r)\to +\infty$ as $r\to R$.
We set $v(x)\decl\varphi(\abs x)$.

Therefore, the function $v$ satisfies the differential equation
$$\diver(\abs{\nabla v}^{p-2} \nabla v)= g(v)\quad\mathrm{on}\quad\Omega_R,$$
where $\Omega_R\decl\{x\ |\ \abs x<R\}$.
On the other hand the function $U\decl-u$ satisfies the inequality 
$$\diver(\abs{\nabla U}^{p-2} \nabla U)\ge g_1(U)\decl g(U)\quad\mathrm{on} \quad \RN.$$ 

Since  $U(x)\le v(x)$ for $\abs x$ close to $R$ we are in the position
to apply the comparison Lemma \ref{lem:comp}. As a consequence,
 $U(x)\le v(x)$ for any $x \in \Omega_R$. In particular $U(0)\le v(0)=a$.
Letting $a\to 0$ it follows that $U(0)\le 0$. Hence $u(0)\ge 0$.

 Next, if $f\ge 0$, then  $u$ is
 a nonnegative  solution of the inequality,
$ -\diver(\abs{\nabla u}^{p-2} \nabla u)\ge0$ on $\RN$.
Hence, by the weak Harnack inequality (see \cite{serrin})
it follows that, either $u\equiv 0$ or $u>0$ on $\RN$.
\ep

\bigskip
The argument for proving Theorem \ref{th:euclmean} is the same of the one used  in the proof of Theorem~\ref{th:main}, so we shall be brief.

\bigskip
\bp{Theorem \ref{th:euclmean}} 
Let $g(t)\decl f(-t)$.
Under the assumptions of Theorem \ref{th:euclmean}, there exists a radial solution 
$v$ of 
$\diver(\displaystyle{\frac{\nabla v}{\sqrt{1+\abs{\nabla v}^2}}})\ge g(v)$
such that $v(0)=a>0$ and $v(r)\rightarrow +\infty$ as $r\rightarrow R$, see \cite{nai-usa97}. By  Lemma \ref{lem:comp} we get $u(0)\ge -a$. Letting $a\rightarrow 0,$
the claim follows.
\ep

\bigskip

\bp{Theorem \ref{th:eqpbd}} Let $u$ be a solution of (\ref{eq:pgen}).
  The function $v:=u-\alpha$ solves the equation
  $$-\Delta_p(v)=f(u-\alpha+\alpha)=f(v+\alpha)=:g(v).$$
  An application of Theorem \ref{th:main} to the last equation implies $v\ge 0$
i.e.  $u\ge \alpha$.

  Next, the function $v:=\beta-u$ solves the equation
  $$-\Delta_p(v)=-f(u)=-f(-v+\beta)=:g(v).$$
  Again, we are in the position to apply Theorem \ref{th:main} which
  yields $v\ge 0$, that is $u\le \beta$. This concludes the proof. 
\ep

\bigskip
\bp{Theorem \ref{th:eqmeanbd}} Step 1.
  We treat first the case when $f$ has at least a zero.  In this case
  let $\alpha$ be the first zero of $f$.
  Let $u$ be a solution of (\ref{dis:mainmean}).
  We set
  \be g(t):=\inf_{s\le t}f(s),\qquad \mathrm{for}\ \  t\in \R. \label{constrg}\ee
  The function $g$ is continuous, non increasing, $g(t)>0$ for $t<\alpha$
  and $g(\alpha)=0$.

  The function $u$ satisfies the inequality 
  $$ -\diver(\frac{\nabla v}{\sqrt{1+\abs{\nabla v}^2}})\ge f(u)\ge g(u), 
      \quad\mathrm{on}\ \ \RN.$$
  By the change of variable $v:=u-\alpha$ and arguing as in the proof of 
  Theorem \ref{th:eqpbd}, we get that $u\ge \alpha.$ The claim is proved.

  Step 2. Next, we assume that $f$ has no zeros. Let $u$ be a solution of
   (\ref{dis:mainmean}). Let $g$ be the function defined in (\ref{constrg}).
  Fix a number $c>0$ and define 
  $$g_c(t):=\min \{ g(t),c-t\}  ,\qquad \mathrm{for}\ \  t\in \R.$$ 
  The function $g_c$  is continuous, non increasing, $g_c(t)>0$ for $t<c$
  and $g_c(c)=0 $.
 Since $u$ satisfies
  $$ -\diver(\frac{\nabla v}{\sqrt{1+\abs{\nabla v}^2}})\ge f(u)\ge g_c(u), 
      \quad\mathrm{on}\ \ \RN,$$
  by step 1 it follows that $u\ge c$. Since the inequality $u\ge c$ holds for any $c$
  we get a contradiction.

  Step 3. Let $u$ be a solution of (\ref{eq:mean1}). 
  By step 1. we have $u\ge \alpha$. In order to prove the estimate $u\le \beta$
 we consider the function $v:= -u$ and argue as in step 1.
  \ep

\bigskip

\bp{Theorem \ref{th:disp}} Arguing by contradiction, let $u$ be a solution of (\ref{dis:main}).
  Fix $\alpha\in\R$ and set $v:=u-\alpha$. The function $v$ solves the inequality
    $-\Delta_p v\ge f(u)=f(v+\alpha)$. Since the function $f(\cdot+\alpha)$ satisfies the hypothesis of
    theorem \ref{th:main} we have $v\ge 0$, that is $u\ge \alpha$. Since the inequality  $u\ge \alpha$ holds for any $\alpha$,
    we obtain $u=+\infty$. This contradiction concludes the proof.\ep
    
    \bigskip

\bp{Theorem \ref{th:meaneq}} Let $l_1:=\lim_{t\to -\infty}f(t)$ and
   $l_2:=\lim_{t\to +\infty}f(t)$.

 We divide the proof into five steps.

 1. We  consider first the case $f>0$ and the inequality (\ref{dis:mean33}).
   In this case $l_1>0$ and Theorem \ref{th:eqmeanbd} yields the thesis.

 2.  Next assume $f<0$ and  let $u$ be a solution of the 
  equation (\ref{eq:mean2}).
  Then $v=-u$ satisfies the equation 
  $\displaystyle -\diver(\frac{\nabla v}{\sqrt{1+\abs{\nabla v}^2}})= g(v)$ where
  $g(t):=-f(-t)$ is a non increasing positive function. 
  This contradicts the step 1.

 Therefore, assume that the function $f$ has at least a zero. 
  Let $u$ be a solution of the  equation (\ref{eq:mean2}).
   From hypotheses on $f$, the set of its zeros, $S$, is an interval
   bounded from one side. 

3.   Assume, first that $S=[\alpha_1,+\infty[$. In this case 
  $l_1>0$ and from Theorem  \ref{th:eqmeanbd} we have that $u\ge \alpha_1$.
  Hence $f(u(x))=0$ for all $x\in\RN$ and so
  $u$ is a function bounded from one side which solves the equation
  \be -\diver(\frac{\nabla u}{\sqrt{1+\abs{\nabla u}^2}})=0,\qquad
  \mathrm{on}\ \ \RN.\label{eq:meanhom}\ee
 This, by Bombieri-De Giorgi-Miranda Theorem (\cite{BDM, DG}) implies that 
 $u\equiv \alpha\ge 0$. 
 
4. Next, assume that $S=]-\infty, \alpha_2]$. By the change of variable
  $v=-u$ and arguing as in the step 3. we obtain that $u\le \alpha_2$.
  Hence $u$ is a function bounded from one side which solves (\ref{eq:meanhom}).
  Therefore $u$ is constant.

5. We analyze the case $S=[\alpha_1, \alpha_2]$. In this case
  $l_1>0$ and $l_2<0$. By Theorem   \ref{th:eqmeanbd} we have that
  $\alpha_1\le u\le \alpha_2$
  that is $u$ takes its value  in $S$. 
  Therefore  $u$ is a  bounded function  which solves (\ref{eq:meanhom}). 
Hence $u$ is constant.
\ep

\section{Some extensions of the main results}
We extend our main results to more general quasilinear operator
in Section \ref{ext}. In Section~\ref{carnot} we extend Theorem \ref{th:main} in 
Carnot group setting. The final Section \ref{porous} deals with a quasilinear inequality
related to the porous medium equation.

\subsection{A class of differential inequalities}\label{ext}
In this section we shall consider inequalities of the type
\be -\diver(A(\abs{\nabla u})\nabla u )\ge f(u)\qquad \mathrm{on}\ \ \RN, \label{dis:gen}
\ee
where we shall assume that
\be \begin{cases} A\in\C(]0,+\infty[),\quad A(t)>0\quad \mathrm{for}\ t>0,\quad \\
  tA(\abs t)\in \C(\R)\cap\Cuno(]0,\infty[)\ \mathrm{and}\   (tA(t))'>0 \quad \mathrm{for}\ t>0.\end{cases}\label{A}
\ee

We shall distinguish two cases accordingly to the asymptotic  behavior of the function $tA(t)$.
Namely,   $\lim_{t\to+\infty}tA(t)<+\infty$ or $\lim_{t\to+\infty}tA(t)=+\infty$.
\bt Let $A$ be  as in (\ref{A}) and such  that $\lim_{t\to+\infty}tA(t)<+\infty$. 
  Let $f:\R\to\R$ be a continuous function satisfying (\ref{cond:f}).
  Let $u$ be a solution of (\ref{dis:gen}), then $u\ge 0$ on $\RN$.
\et
We observe that Theorem \ref{th:euclmean} is a particular case of the above theorem.

In the case  $\lim_{t\to+\infty}tA(t)=+\infty,$  then accordingly with \cite{NS1,NS2}, we can construct the following function $G$ and $H$.
Let $G$ be defined as
$$G(t)\decl t^2 A(t)-\int_0^t s A(s) ds,\quad  t\ge 0.$$
The function $G$ is continuous, strictly increasing,  $G(0)=0$ and $G(+\infty)=+\infty$,
see \cite{NS1,NS2,nai-usa97}. Let $H$ be its inverse: the function $H$ is increasing and 
$H(+\infty)=+\infty$.

\bt\label{th:maingen}  Let $A$ be as in (\ref{A}) and such  that 
  $\lim_{t\to+\infty}tA(t)=+\infty$. 
  Let $f:\R\to\R$ be a continuous function satisfying (\ref{cond:f}) and 
 \be \int_{-\infty}^{-1} \frac{1}{H\left(\int_t^{-1} f(s)\ ds\right)} \ dt< +\infty. 
      \label{dis:naiusagen}\ee
  Let $u$ be a solution of (\ref{dis:gen}), then $u\ge 0$ on $\RN$.
\et

If $A(t)=t^{p-2}$, then  $H(t)=\left(\frac{p}{P-1}\right)^{1/p} t^{1/p} $ 
and the above theorem is indeed Theorem \ref{th:main}.

We leave the proof of the above results to the interested reader since it is  based on the same idea
already discussed above by taking into account the nonexistence Theorem 1 and Theorem 2 in \cite{nai-usa97}.

Condition (\ref{dis:naiusagen}) is sharp in the following sense. 
If the integral in  (\ref{dis:naiusagen}) diverges,  then (\ref{dis:gen}) admits a negative solution.
Indeed from  Theorem 3 in \cite{nai-usa97} it follows  that  equation 
(\ref{eq:radial}) with $g(t)=f(-t)$ has a positive solution. This implies the claim.
\begin{example} Let $A(t)\decl\frac{\ln(1+t)}{t}$ and $f(s)\ge c\abs s^q$ for $s<0$.
  We claim that if $q>0$ the solutions of
  \be -\diver\left( {\ln(1+\abs{\nabla u})}\frac{\nabla u}{\abs {\nabla u}}\right)\ge f(u)
  \qquad \mathrm{on}\ \ \RN
    \label{dis:example}\ee
  are nonnegative.
  We are in the position to apply
  Theorem \ref{th:maingen}. In this case $G(t)= t-\ln(1+t)$.  
  In order to prove the claim it is enough to show that the function $H(T)$ behaves at infinity as $T$.
  Since $H(T)$ is the solution of $t-\ln(1+t)=T$ by the change of variable
  $z=(1+T)^{-1}$ and   $x=(1+t)^{-1}$ the equation becomes
  $z=\frac{x}{1+x\ln x}$ and we have to study the zero of the function $x(z)$ defined implicitly.
  It is easy to recognize that $x(z)=z+o(z)$, and this  implies the claim.

  Using the same argument as above  one can prove the following. We omit the details.
\end{example}
\bt Let $f:\R\to\R$ be a continuous function such that : there exists $c>0$ and  $q>1$ such that
   $f(s)\ge\, c\,(\ln(1+\abs s))^q$ for $s<0$.  Then  the solutions of (\ref{dis:example}) are nonnegative.
\et

The analogous results of Theorem \ref{th:eqpbd}, \ref{th:eqmeanbd},
\ref{th:eqp}, \ref{th:disp}, \ref{th:meaneq}
are the following.

\bt\label{th:eqpbdA}  Let $A$ be  as in (\ref{A}) and such  that 
  $\lim_{t\to+\infty}tA(t)=+\infty$. 
Let $f\colon\R\to\R$ be a continuous function such that there exists $\alpha,\beta\in\R$, $\alpha\le\beta$ such that
\be f_{]-\infty,\alpha[}\quad \mathrm{ is\ positive\ and\ non\ increasing,}\quad 
f_{]\beta,+\infty[}\quad \mathrm{ is\ negative\ and\ non\ increasing}, 
    \label{cond:foddA}\ee
and
\be \int_{-\infty}^{\alpha} \frac{1}{H\left(\int_t^{\alpha} f(s)\ ds\right)} \ dt< +\infty , \quad
   \int_\beta^{\infty} \frac{1}{H\left(\int_\beta^{t} -f(s)\ ds\right)} \ dt< +\infty. \label{dis:naiusa00G}\ee

  Let $u$ be a solution of 
  \be -\diver(A(\abs{\nabla u})\nabla u )= f(u)\qquad \mathrm{on}\ \ \RN, \label{eq:gen}\ee
  then $u$ is bounded and  $\alpha\le{u(x)}\le \beta$  for any $x\in \RN$.
\et

\bt\label{th:eqmeanbdA}  Let $A$ be as in (\ref{A}) and such  that $\lim_{t\to+\infty}tA(t)<+\infty$. 
  Let $f\colon\R\to\R$ be a continuous function such that 
  $$\liminf_{t\to-\infty} f(t)>0.$$
  If $u$ is a solution of (\ref{dis:gen}), then
  $f$ has at least a zero, and set $\alpha:=$ the first zero of $f$
  (that is $\alpha:=\min S$ where $S:= f^{-1}(0)$) we have $u\ge \alpha$.
  In particular if $f>0$ the (\ref{dis:gen}) has no solution.

  Moreover if   $$\limsup_{t\to+\infty} f(t)<0$$
  and $u$ solves (\ref{eq:gen})
  then $u$ is bounded and  $\alpha\le{u(x)}\le \beta$  for any $x\in \RN,$
  where $\beta:=$ last zero of $f$ (that is $\beta:=\max S$).
\et

\bt\label{th:eqpA}  Let $A$ be  as in (\ref{A}) and such  that 
  $\lim_{t\to+\infty}tA(t)=+\infty$. 
Let $f\colon\R\to\R$ be a non increasing continuous function such that 
\be f(t)>0\quad \mathrm{ if}\quad  t<0, \ \  \mathrm{and}\  \  f(t)<0\quad \mathrm{ if}\quad  t>0, \label{cond:oddA}\ee
and (\ref{dis:naiusa00G}) holds with $\alpha=-1$ and $\beta=1$.

  If $u$ be a solution of (\ref{eq:gen}),   then $u\equiv 0$ on $\RN$.
\et

\bt\label{th:dispA}  Let $A$ be  as in (\ref{A}) and such  that 
  $\lim_{t\to+\infty}tA(t)=+\infty$. 
  Let $f\colon\R\to\R$ be a positive, non increasing, continuous function 
satisfying (\ref{dis:naiusagen}).
Then the inequality (\ref{dis:gen}) has no solutions.
\et

\bt\label{th:meaneqA}  Let $A$ be as in (\ref{A}) and such  that 
  $\lim_{t\to+\infty}tA(t)<+\infty$. 
 Let $f\colon\R\to\R$ be a non increasing,
  continuous function. If $u$ is a solution of 
  (\ref{eq:gen})  then
  $u$ solves the homogeneous problem
   \be -\diver(A(\abs{\nabla u})\nabla u )= 0\qquad \mathrm{on}\ \ \RN, \label{eq:genhom}\ee
  and for all $x\in \RN$, $u(x)\in S$ where $S:=f^{-1}(0)$.
In particular if $f(t)\neq 0$ for any $t$, then (\ref{eq:gen}) has no solutions.

In addition, if $f$ is supposed to be positive then the inequality
(\ref{dis:gen}) has no solutions.
\et
\boss
Very recently, by using a completely different technique,
 a strong generalization of Theorem \ref{th:meaneqA} 
has been obtained by James Serrin
\cite[Theorem 3.23]{serrin1}. 
We wish to thank  Alberto Farina for pointing out this information.
\eoss

For   recent nonexistence results related to anticoercive problems,
 we refer the interested reader to the forthcoming paper \cite{analyse}. 

\subsection{Inequalities on Carnot Groups}\label{carnot}
Let $\RN\equiv\G$ be a Carnot group and let $\grl$ be the horizontal gradient on $G$ and $Q>1$ the homogeneous dimension (see Appendix and  \cite{bonfi} for details on these structures).
Let $\Gamma_p$ be the fundamental solution of the quasilinear operator
$\Delta_{L,p}u=\diverl(\abs{\grl u}^{p-2}\grl u)$ at the origin.
Set
$$ N_p\decl \begin{cases}
               \Gamma_p^{\frac{p-1}{p-Q}} &  p>1,p\neq Q \cr
               \exp (-\Gamma_p)     & p=Q\end{cases}$$
It is known that $N_p$ is a homogeneous norm on 
$\G$. From \cite{capogna, cap-dan-gar96}, it is known that  $N_p$ is 
 H\"older continuous.
In what follows we shall assume that $N_p$ is smooth. This assumption is satisfied for example for ``Heisenberg type'' groups. See  \cite{cap-dan-gar96}.

With the above notation, we have that if $\zeta\colon\R\to\R$ is a smooth function,
then the \emph{radial} function $v\decl\zeta\circ N_p\colon\G\to \R$ satisfies
$$\Delta_{G,p}v:=\diverl(\abs{\grl v}^{p-2} \grl v)=
(p-1)\psi^p\abs{\zeta'}^{p-2}\left(\zeta''(r)+\frac{Q-1}{p-1}\frac{\zeta'(r)}{r} 
   \right)_{r=N_P},$$
where $\psi\decl \abs{\grl N_p}$ is a bounded function, see \cite{dam05}.
Hence we can apply the same arguments used in the preceding section 
obtaining an analog of Theorem \ref{th:main} in this more general setting.

\bt\label{th:mainG} Let $p>1$. Let $f\colon\R\to\R$ be a continuous function  
 satisfying  (\ref{cond:f}) and (\ref{dis:naiusa}).
Let $u$ be a solution of 
  \be -\diverl(\abs{\grl u}^{p-2} \grl u)\ge f(u)\quad\mathrm{on} \quad \RN,\label{dis:mainG} \ee
  then $u\ge 0$ on $\RN$. Moreover if $f(t)\ge 0$ for $t\ge0$ then, either $u\equiv 0$ or $u>0$ on $\RN$.
\et

Immediate consequences of  Theorem \ref{th:mainG} 
are the following Liouville type theorems.
\bt\label{th:eqpG}   Let $p>1$. Let $f\colon\R\to\R$ be a continuous function such that there exists $\alpha,\beta\in\R$, $\alpha\le\beta$ such that
(\ref{cond:fodd}) and (\ref{dis:naiusa00}) hold.
  If $u$ is a solution of 
  \be -\diverl(\abs{\grl u}^{p-2} \grl u)=  f(u)\quad\mathrm{on} \quad \RN, 
    \label{eq:plapgG}\ee
  then $u$ is bounded and  $\alpha\le{u(x)}\le \beta$  for any $x\in \RN$.
In particular, if $q>p-1$ and $u$ is a solution of 
  \be
   \Delta_{G,p} u= \abs{u}^{q-1}u\quad\mathrm{on} \quad \RN, \label{breziscarnot}
   \ee
  then $u\equiv 0$ on $\RN$.
\et
\boss
The conclusion related to inequality (\ref{dis:mainG})
in the Euclidean setting,  $p=2,$  and $f(u) = \abs{u}^{q-1}u$ has been obtained by Brezis \cite{brezis}
by using a variant of Kato's inequality and .
\eoss

\bt\label{th:dispG} Let $p>1$. Let $f\colon\R\to\R$ be a positive, non increasing, continuous function 
satisfying (\ref{dis:naiusa}).
Then the inequality (\ref{dis:mainG}) has no solutions.
\et 

\bc\label{cor:liouvG} Let $p\ge Q>1$ and $f:\R\to[0,+\infty[$ be  a continuous function satisfying (\ref{cond:f}) and
(\ref{dis:naiusa}). If $u$ is a solution of  (\ref{dis:main}) then $u$ is constant on $\RN$.
More precisely $u\equiv \alpha\ge 0$ and $f(\alpha)=0$.

Moreover  if $f(t)>0$ for $t\ge0$, then  (\ref{dis:main}) has no solutions.
\ec

\bp{Theorem \ref{th:mainG}} Let $u$ be a solution of (\ref{dis:main}). 
Since the inequality is invariant under translations, it is sufficient to prove that $u(0)\ge 0$.

  Let $C>0$ be a constant such that $\psi^p\le C$. Set $g(t)\decl f(-t)/C$. The function $g$ satisfies the assumptions of Theorem \ref{th:blowup}.
Let $D=Q>1$ be the homogeneous dimension. 
Let $a>0$ and let $\varphi$ be a solution of (\ref{eq:radial}) such that 
$\varphi(r)\to +\infty$ as $r\to R$.
We set $v(x)\decl\varphi(N_p(x))$.
By computation we have, 
\bern \diverl(\abs{\grl v}^{p-2} \grl v)=(p-1)\psi^p\abs{v'}^{p-2}\left(v''(r)+\frac{Q-1}{p-1}\frac{v'(r) }{r}
   \right)_{r=N_p}\\
=\psi^pN_p^{1-Q} \left(r^{Q-1} \abs{v'(r)}^{p-2}v'(r)\right)'_{r=N_p}.\eern
Therefore, the function $v$ satisfies the differential equation
$$\diverl(\abs{\grl v}^{p-2} \grl v)= g_2(x,v)\decl \psi^p g(v)\le C g(v)$$
on $\Omega_R\decl\{x\ |\ N_p(x)<R\}$.
On the other hand the function $U\decl-u$ satisfies the inequality 
$$\diverl(\abs{\grl U}^{p-2} \grl U)\ge f(-U)=Cg(U)=:g_1(U)\quad\mathrm{on} \quad \RN.$$ 

Since $g_1\ge g_2$ and $U(x)\le v(x)$ for $N_p(x)$ close to $R$ we are in the position
to apply the comparison Lemma \ref{lem:comp}. As a consequence,
 $U(x)\le v(x)$ for any $x \in \Omega_R$. In particular $U(0)\le v(0)=a$.
Letting $a\to 0$ it follows that $U(0)\le 0$. Hence $u(0)\ge 0$.

 Next, if $f\ge 0$, then  $u$ is
  nonnegative and solves the inequality
$ -\diverl(\abs{\grl u}^{p-2} \grl u)\ge0$ on $\RN$.
Hence, by the weak Harnack inequality (see \cite{capogna})
it follows that, either $u\equiv 0$ or $u>0$ on $\RN$.
\ep
\medskip

The proofs of the above  Liouville theorems  are very similar to those given in 
Section \ref{sec:main}, so we omit them. 
The proof of Corollary \ref{cor:liouvG} relies the fact that the only
nonnegative functions $u$ such that 
$-\Delta_{G,p} u\ge0$ on $\RN$ with $Q\le p$
are the constants, see \cite{dam09}.

\subsection{A  porous medium type inequality}\label{porous}
We end this paper by pointing out the following  slight modification of Theorem  \ref{th:main}.
\bt Let $\gamma\ge 1$.  Let $f:\R\to\R$ be a continuous function satisfying (\ref{cond:f}) and
 \be \int_{-\infty}^{-1} \abs{t}^{\gamma-1}  {\left(\int_t^{-1} f(s)\abs{s}^{\gamma-1}\ ds\right)^{-\frac12}} \ dt< +\infty. \label{dis:fast}\ee
   Let $u$ be a solution of
 $$ -\Delta(\abs u^{\gamma-1}u) \ge f(u)\qquad \mathrm{on}\ \ \RN,$$ 
 then $u\ge 0$ on $\RN$.

  In particular if $u$ is a solution and  $f(t)\ge C\abs t^q$ for $t<0$ with $q> \gamma$, then 
   either $u\equiv 0$ or $u>0$ on $\RN$.
\et 

\section*{Appendix}
We quote some  facts on Carnot groups  and refer the interested reader to  
\cite{bonfi,fol-ste82} 
  for  more detailed information on this subject.

A Carnot group is a connected, simply connected, nilpotent Lie
group $\G$ of dimension $N$ with graded Lie algebra ${\cal
G}=V_1\oplus \dots \oplus V_r$ such that $[V_1,V_i]=V_{i+1}$ for
$i=1\dots r-1$ and $[V_1,V_r]=0$. Such an integer $r$ is called the
\emph{step} of the group.
 We set $l=n_1=\dim V_1$, $n_2=\dim V_2,\dots,n_r=\dim V_r$.
A  Carnot group $\G$ of dimension $N$ can be identified, up to an
isomorphism, with the structure of a \emph{homogeneous Carnot
Group} $(\RN,\circ,\delta_R)$ defined as follows; we
identify $\G$ with $\RN$ endowed with a Lie group law $\circ$. We
consider $\RN$ split in $r$ subspaces
$\RN=\R^{n_1}\times\R^{n_2}\times\cdots\times\R^{n_r}$ with
$n_1+n_2+\cdots+n_r=N$ and $\xi=(\xi^{(1)},\dots,\xi^{(r)})$ with
$\xi^{(i)}\in\R^{n_i}$. 
We shall assume that for any $R>0$ the dilation
$\delta_R(\xi)=(R\xi^{(1)},R^2 \xi^{(2)},\dots,R^r \xi^{(r)})$ 
is a Lie group automorphism.
The Lie algebra of
left-invariant vector fields on $(\RN,\circ)$ is $\cal G$. For
$i=1,\dots,n_1=l $ let $X_i$ be the unique vector field in $\cal
G$ that coincides with $\partial/\partial\xi^{(1)}_i$ at the
origin. We require that the Lie algebra generated by
$X_1,\dots,X_{l}$ is the whole $\cal G$.

We denote with $\grl$ the vector field $\grl\decl(X_1,\dots,X_l)^T$
and we call it \emph{horizontal vector field}.
Moreover, the vector
fields $X_1,\dots,X_{l}$ are homogeneous of degree 1 with respect
to $\delta_R$ and in this case
$Q= \sum_{i=1}^r i\,n_i=  \sum_{i=1}^r i\,\mathrm{dim}V_i$ is called the 
\emph{homogeneous dimension} of $\G$.
The \emph{canonical sub-Laplacian} on $\G$ is the
second order differential operator defined by
$\Delta_G=\sum_{i=1}^{l} X_i^2$ and for $p>1$ the
$p$-sub-Laplacian operator is
$\sum_{i=1}^{l} X_i(\abs{\grl u}^{p-2}X_iu)$.
Since $X_1,\dots,X_{l}$
generate the whole $\cal G$, the sub-Laplacian $\Delta_G$ satisfies the
H\"ormander hypoellipticity  condition.

A nonnegative continuous function $N:\RN\to\R_+$ is called a 
\emph{homogeneous norm} on {\doppio G}, if 
$N(\xi^{-1})=N(\xi)$, $N(\xi)=0$ if and only if $\xi=0$, and it is
homogeneous of degree 1 with respect to $\delta_R$ (i.e.
$N(\delta_R(\xi))=R N(\xi)$).
A homogeneous norm $N$ defines on $\G$ a \emph{pseudo-distance} defined as
$d(\xi,\eta)\decl N(\xi^{-1}\eta)$, which
in general is not a distance.
If $N$ and $\tilde N$ are two homogeneous norms, then they are equivalent,
that is, there exists a constant
$C>0$ such that $C^{-1}N(\xi)\le \tilde N(\xi)\le CN(\xi)$.
Let $N$ be a homogeneous norm, then there exists a constant
$C>0$ such that $C^{-1}\abs\xi\le N(\xi)\le C\abs\xi^{1/r}$,
for $N(\xi)\le1$.
An example of homogeneous norm is 
$  N(\xi)\decl\left(\sum_{i=1}^r\abs{\xi_i}^{2r!/i}\right)^{1/2r!}.$

Notice that if $N$ is a homogeneous norm differentiable a.e., 
then $\abs{\grl N}$ is homogeneous of degree 0 with respect to 
$\delta_R$; hence $\abs{\grl N}$ is bounded.

\medskip

Special examples of Carnot groups are the
  Euclidean spaces $\R^Q$.
  Moreover, if $Q\le 3$ then any Carnot group is the ordinary Euclidean
  space $\R^Q$.

 The most simple nontrivial example of a Carnot group
  is the Heisenberg group $\hei^1=\R^3$.
  For an integer $n\ge1$, the Heisenberg group $\hei^n$ is defined as follows:
  let $\xi=(\xi^{(1)},\xi^{(2)})$ with
  $\xi^{(1)}\decl(x_1,\dots,x_n,y_1,\dots,y_n)$ and $\xi^{(2)}\decl t$.
  We endow $\R^{2n+1}$ with  the group law
$\hat\xi\circ\tilde\xi\decl(\hat x+\tilde x,\hat y+\tilde y,\hat t+ \tilde t+2\sum_{i=1}^n(\tilde x_i\hat y_i-\hat x_i \tilde y_i)).$
We consider the vector fields
\[X_i\decl\frac{\partial}{\partial x_i}+2y_i\frac{\partial}{\partial t},\
        Y_i\decl\frac{\partial}{\partial y_i}-2x_i\frac{\partial}{\partial t},
        \qquad\mathrm{for\ } i=1,\dots,n, \]
and the associated  Heisenberg gradient
$ \grh\decl (X_1,\dots,X_n,Y_1,\dots,Y_n)^T$.
The Kohn Laplacian $\lh$ is then the operator defined by
$\lh\decl\sum_{i=1}^nX_i^2+Y_i^2.$
The family of dilations is given by
$\delta_R(\xi)\decl (R x,R y,R^2 t)$ with homogeneous dimension
$Q=2n+2$.
In ${\hei^n}$ a canonical homogeneous norm is defined as
$\abs{\xi}_H\decl \left(\left(\sum_{i=1}^n x_i^2+y_i^2\right)^2+t^2\right)^{1/4}.$

\def\cprime{$'$} \def\cprime{$'$} \def\cprime{$'$} \def\cprime{$'$}
  \def\cprime{$'$} \def\cprime{$'$}
  
  \bigskip
  
\noindent{\bf Acknowledgements} We wish to thank Alberto Farina and James Serrin  for  useful discussions and for attracting our attention to references
\cite{brezis}, \cite{FA}, \cite{FAS}, \cite{serrin1}.

\end{document}